\documentclass[12pt]{amsart}

\usepackage{amssymb}
\usepackage{graphicx}
\usepackage{psfrag}
  
\theoremstyle{plain}
\newtheorem{theorem}{Theorem}[section]
\newtheorem{lemma}[theorem]{Lemma}
\newtheorem{proposition}[theorem]{Proposition}
\newtheorem{corollary}[theorem]{Corollary}

\theoremstyle{definition}
\newtheorem{definition}[theorem]{Definition}
\newtheorem{algorithm}[theorem]{Algorithm}

\theoremstyle{remark}
\newtheorem{remark}[theorem]{Remark}

\renewenvironment{proof}[1]{\vspace*{.1in}\noindent{\bf
Proof{#1}. \/}}{\qed\vspace{3ex}} 
		     
\newcommand{\bbZ}{{\mathbb Z}}
\newcommand{\bbH}{{\mathbb H}}

\newcommand{\bH}{{\mathbf H}}
\newcommand{\blambda}{{\boldsymbol{\lambda}}}

\newcommand{\cF}{{\mathcal F}}

\renewcommand{\tilde}{\widetilde}
\newcommand{\stab}{\operatorname{Stab}}

\begin{document}
\title{Small covers of the dodecahedron and the $120$-cell}

\def \makeauthor{
\author{A.~Garrison and R.~Scott}
\address{Department of Mathematics\\
Santa Clara University\\ 
Santa Clara, CA 95053}
\email{rscott@math.scu.edu}
} 
\makeauthor

\date{July, 2001}

\begin{abstract}
Let $P$ be the right-angled hyperbolic dodecahedron or $120$-cell, and
let $W$ be the group generated by reflections across codimension-one
faces of $P$.  We prove that if $\Gamma\subset W$ is a torsion free 
subgroup of minimal index, then the corresponding hyperbolic manifold
$\bbH^n/\Gamma$ is determined up to homeomorphism by $\Gamma$ modulo
the symmetry group of $P$.
\end{abstract}

\maketitle

\section{Introduction}

Let $P$ be a right-angled bounded convex polytope in hyperbolic space
$\bbH^n$, and let $W$ be the group generated by reflections across 
codimension-one faces.  For any torsion-free subgroup $\Gamma\subset
W$ of finite index, the quotient $\bbH^n/\Gamma$ is a closed
hyperbolic manifold which is an orbifold cover of $P$.  We call this
manifold a {\em small cover of $P$} (as in \cite{DJ}) if the index of
$\Gamma$ is minimal.  Examples of small covers include the first closed
hyperbolic $3$-manifold to appear in the literature \cite{Lo} as well 
as its generalizations \cite{Ve}.

The point of this note is to show that if $P$ is regular then any
small cover of $P$ is uniquely determined up to isometry
(and up to homeomorphism when $n\geq 3$) by $\Gamma$ (modulo 
symmetries of $P$).  In fact, there are only two right-angled regular
hyperbolic polytopes with dimension $\geq 3$, the dodecahedron and the
$120$-cell.  We conclude the
paper by showing that up to homeomorphism there are exactly $25$ small
covers of the dodecahedron, and that there is a unique small cover of
the $120$-cell with minimal complexity (in the sense of
Section~\ref{s:examples}, below). 

There are other hyperbolic manifolds based on the dodecahedron
\cite{SW} and the $120$-cell \cite{D} appearing in the literature, but
these are obtained by identifying faces of a single copy of the
polytope and require that it have dihedrals angles of $2\pi/5$.
Analogous constructions fail for right-angled realizations of these
polytopes, and small covers provide a natural alternative. 

\section{Definitions and preliminary facts}

Let $P$ be an $n$-dimensional right-angled convex polytope in
$\bbH^n$, and let $\cF$ denote the set of facets (i.e.,
codimension-one faces) of $P$.  For each $F\in \cF$, we let $s_F$
denote the reflection across $F$, and we let $W$ be the group
generated by $\{s_F\;|\;F\in\cF\}$.  Defining relations for $W$ are
$(s_F)^2=1$ for all $F$ and $s_Fs_{F'}=s_{F'}s_F$ whenever $F\cap F'\neq
\emptyset$.  Following \cite{DJ}, we call an
epimorphism $\lambda:W\rightarrow(\bbZ_2)^n$ a {\em characteristic
function} if whenever $F_1,\ldots,F_n$ are facets that all meet at a
vertex, the images $\lambda(s_{F_1}),\ldots,\lambda(s_{F_n})$ form a
$\bbZ_2$-basis.    

\begin{proposition}
If $\Gamma$ is a torsion free subgroup of $W$, then its index is $\geq
2^n$.  If the index is equal to $2^n$, then $\Gamma$ is normal and is
the kernel of a characteristic function
$\lambda:W\rightarrow(\bbZ_2)^n$.  
\end{proposition}

\begin{proof}{}
We consider the reflection tiling of $\bbH^n$ corresponding to $P$ and
$W$.  If $C$ is a codimension-$k$ cell in this tiling, then there are
exactly $2^k$ maximal cells containing $C$, and the
stabilizer $\stab(C)$ is isomorphic to $(\bbZ_2)^k$.  If $\Gamma$ 
is torsion-free, then the natural map $\stab(C)\rightarrow W/\Gamma$
must be injective; hence, the index of $\Gamma$ is at least $2^n$.     

If the index is equal to $2^n$, then for any $0$-cell $v$, the map
$\phi:\stab(v)\rightarrow W/\Gamma$ is a bijection.
We let $\lambda$ be the composition 
\[W\longrightarrow
W/\Gamma\stackrel{\phi^{-1}}{\longrightarrow}\stab(v)\cong(\bbZ_2)^n.\]   
It is not {\em a priori} a homomorphism, 
but since the defining relations for $W$ also hold for
the images $\{\lambda(s_F)\;|\; F\in\cF\}$, there is an induced
epimorphism $\mu:W\rightarrow\stab(v)$ defined by
$\mu(s_F)=\lambda(s_F)$, for all $F\in\cF$.  It is clear that if
$\lambda$ and $\mu$ agree on all but one element in the stabilizer
of a cell, then they must also agree on the entire stabilizer.  Using
this fact and the fact that $\lambda$ and $\mu$ agree on the
generators of $W$, it follows that $\lambda=\mu$.  In particular,
$\lambda$ is an epimorphism and $\Gamma$ is its kernel.  That
$\lambda$ is a characteristic function follows from the fact that 
$\stab(v)\rightarrow W/\Gamma$ must be a bijection {\em for every $v$}.
\end{proof}

\begin{definition}
Let $\lambda:W\rightarrow(\bbZ_2)^n$ be a characteristic function, and
let $\Gamma_{\lambda}=\ker(\lambda)$.  The {\em
small cover of $P$ associated to $\lambda$} is the closed
hyperbolic manifold $M_{\lambda}=\bbH^n/\Gamma_{\lambda}$.
\end{definition}

Small covers are functorial in the following sense.

\begin{proposition}\label{prop:subface}
Let $M_{\lambda}$ be a small cover of $P$.  For any face $F$ of $P$,
let $\bbH_F\subset\bbH^n$ be the hyperbolic subspace spanned by 
$F$ and let $M_F$ be the image of $\bbH_F$ in $M_{\lambda}$ (that is,
$M_F=\bbH_F/\left(\Gamma_{\lambda}\cap \stab(\bbH_F)\right)$).  Then
$M_F$ is a totally geodesic submanifold which is itself a small cover of
the face $F$. 
\end{proposition}

\begin{proof}{}
The stabilizer $\stab(F)\subset W$ fixes
$\bbH_F$ pointwise and is mapped isomorphically by $\lambda$ onto a
codimension-$k$ subspace $V_F\subset(\bbZ_2)^{n}$ (where
$k=\dim(F)$).  Each facet of $F$ can 
be expressed uniquely as the intersection $F'\cap F$ where $F'$ is a
facet of $P$ that is orthogonal to $F$, and we let $W_F$ denote the
subgroup of $W$ generated by the corresponding reflections $s_{F'}$.   
The characteristic function $\lambda$ induces a characteristic
function $\overline{\lambda}:W_F\rightarrow
(\bbZ_2)^n/V_F$ for the polytope $F$ with kernel
$\Gamma_{\lambda}\cap\stab(\bbH_F)$.  It follows that
$M_F=M_{\overline{\lambda}}$.  
\end{proof}  

Let $A$ ($=A(P)$) denote the symmetry group of $P$, and let
$G=GL_n(\bbZ_2)$.  We say that two small covers $M_{\lambda}$ and
$M_{\mu}$ are {\em equivalent} if $\lambda=g\circ\mu\circ a$ for some
$a\in A$ and $g\in G$.  Any small cover $M_{\lambda}$ has a natural 
(totally geodesic) cell decomposition induced by the tiling of
$\bbH^n$.

\begin{proposition}\label{prop:equivalence}
Two small covers $M_{\lambda}$ and $M_{\mu}$ are equivalent if and
only if they are isomorphic as cell complexes.  
\end{proposition}

\begin{proof}{}
It is clear that equivalent small covers are isomorphic as cell
complexes.  Conversely, suppose $\phi:M_{\lambda}\rightarrow M_{\mu}$
is an isomorphism of cell complexes.  Then the lift
$\tilde{\phi}:\bbH^n\rightarrow\bbH^n$ is an automorphism of the
tiling of $\bbH^n$ such that
$\Gamma_{\lambda}=\tilde{\phi}\Gamma_{\mu}\tilde{\phi}^{-1}$. 
Since the automorphism group of the tiling is a semi-direct product of
$W$ and $A$, there exists an $a\in A$ such that
$\Gamma_{\lambda}=w\Gamma_{\mu\circ a}w^{-1}$ for some element $w\in
W$.  It follows that $\Gamma_{\lambda}=\Gamma_{\mu\circ a}$ (they are
both normal), hence $\lambda=g\circ\mu\circ a$ for some $g\in G$. 
\end{proof}

\section{The main theorem}

Let $P$ be a right-angled hyperbolic polytope as above.  $P$ is {\em
regular} if its symmetry group acts transitively on the simplices in
its barycentric subdivision.  Any such simplex must be a hyperbolic   
Coxeter simplex, and by examining the standard list of these simplices
(Bourbaki \cite{B}, p. 133) one can show that $P$ must be either a
polygon with $>4$ sides, the dodecahedron, or the ($4$-dimensional)
$120$-cell.  An important property of these regular polytopes is the
following:

\begin{lemma}\label{lem:connectingpath}
Let $P$ be a regular right-angled hyperbolic polytope, and let $L$ be
the length of an edge in $P$.  Then any path in $P$ that connects
non-adjacent facets has length $\geq L$ with equality holding only if
the connecting path is an edge of $P$.
\end{lemma}  

\begin{proof}{}
Let $\alpha$ be a connecting path with endpoints lying on the
(disjoint) facets $F$ and $F'$.  Because $P$ is regular, it has a 
``cubical'' decomposition obtained by cutting each edge with the orthogonal
hyperplane through its midpoint (see Figure~\ref{fig:Pminpath.fig}).
\begin{figure}[ht]
\begin{center}
\psfrag{H}{$H_E$}
\psfrag{D}{$D$}
\psfrag{B}{$B$}
\psfrag{a}{$\alpha$}
\psfrag{F}{$F$}
\psfrag{F'}{$F'$}
\psfrag{E}{$E$}
\includegraphics[scale = .3]{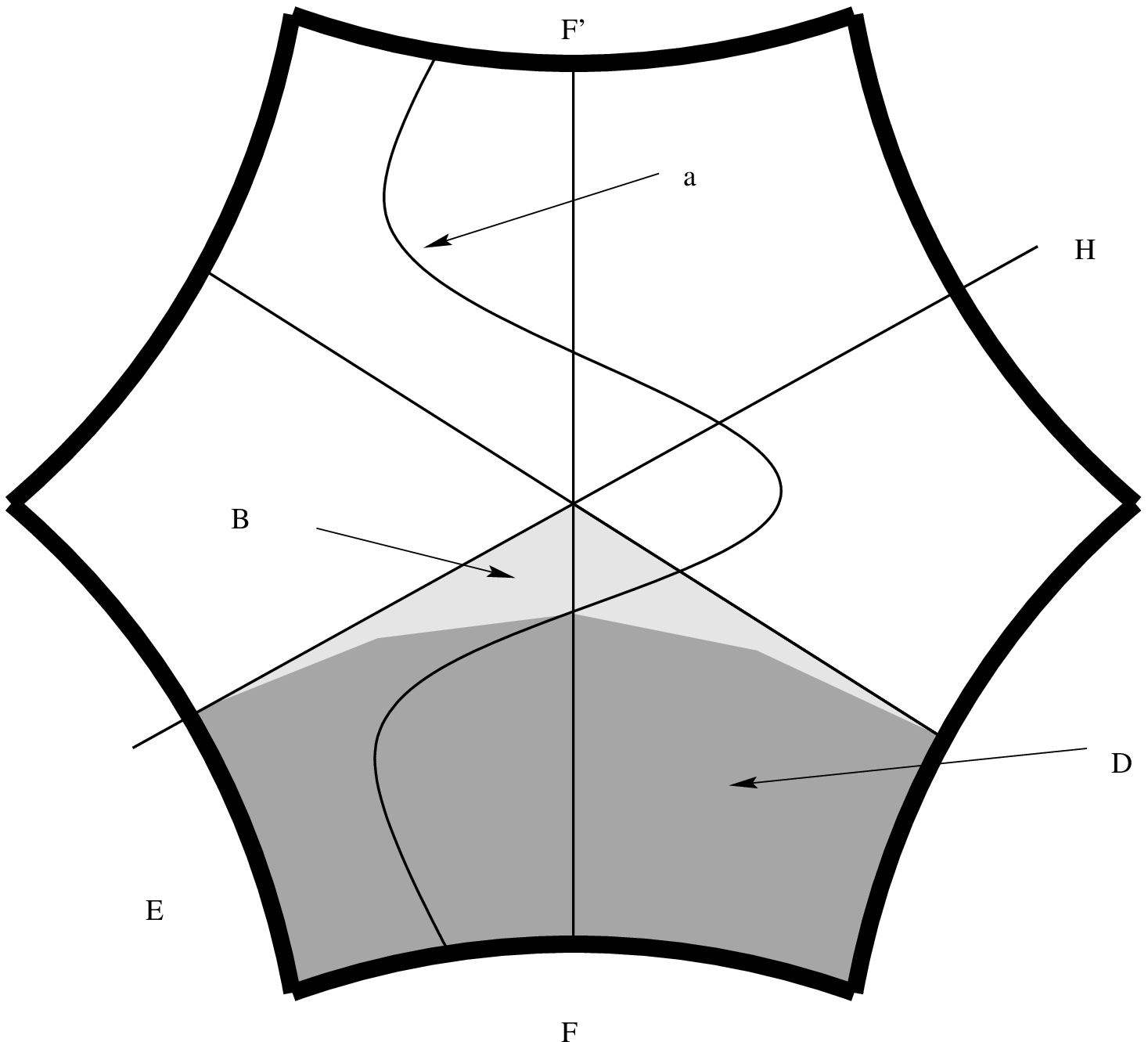}
\caption{\label{fig:Pminpath.fig}}
\end{center}
\end{figure}
We define two open neighborhoods $B$ and $D$ of the facet $F$ as 
follows.  $B$ is the union of the interiors of all cubes that
intersect $F$, and $D$ is the set of all points of $P$ whose distance
to $F$ is less than $L/2$.  We claim that $D\subset B$.  To see this,
it suffices to note that if $E$ is an edge meeting $F$ othogonally at
a vertex and $H_E$ is 
the hyperplane bisecting $E$, then $D$ lies entirely on one
side of $H_E$ (since the edge $E$, being a common perpendicular to $H_E$
and $F$, realizes the shortest distance between the corresponding
hyperplanes).  Similarly, if $B'$ and $D'$ are the corresponding
neighborhoods of $F'$, we have $D'\subset B'$.  Since $F$ and $F'$ are
not adjacent, $B\cap B'=\emptyset$ and, hence, $D\cap D'=\emptyset$.
This means
\[l(\alpha)\geq l(\alpha\cap D)+l(\alpha\cap D')\geq L/2+L/2=L,\] 
and since the edge $E$ is the {\em unique} path of minimal
length joining $F$ to $H_E$, the equality $l(\alpha)=L$ is only
possible if $\alpha$ is an edge of $P$.  
\end{proof}

\begin{remark} 
To see that the lemma can fail without the regularity
hypothesis, let $P$ be any right-angled $n$-dimensional polytope and
let $E\subset P$ be a minimal length edge.  By gluing $2^{n-1}$ copies
of $P$ together around the edge $E$, one obtains a right-angled
polytope $Q$ with $E$ now being an interior connecting path that is as
short as any edge of $Q$. 
\end{remark}

We are now in a position to prove the main theorem.  The argument is
based on the fact that an isometry must preserve the set of minimal
length closed geodesics.  To simplify the exposition, we call a closed
geodesic in a small cover an {\em edge loop} if it is of the form
$M_E$ (Proposition~\ref{prop:subface}) for some edge $E\subset P$. 

\begin{theorem}\label{thm:main}
Let $P$ be an $n$-dimensional right-angled regular hyperbolic
polytope.  Then two small covers of $P$ are isometric if and only 
if they are equivalent. 
\end{theorem}

\begin{proof}{}
First we show that if $M$ is any small cover of $P$, then any closed
geodesic of minimal length must be an edge loop.  Let
$\alpha$ be a closed geodesic in $M$, and let $F$ be a codimension-one
cell in $M$ that intersects $\alpha$ transversely.  Lifting $\alpha$
to the universal cover $\bbH^n$, we obtain a geodesic segment
$\tilde{\alpha}$ that connects two lifts $\tilde{F}$ and $\tilde{F}'$
of the cell $F$ (see Figure~\ref{fig:ucover.fig}).  
\begin{figure}[ht]
\begin{center}
\psfrag{H}{$H$}
\psfrag{H'}{$H'$}
\psfrag{F}{$\tilde{F}$}
\psfrag{F'}{$\tilde{F}'$}
\psfrag{a}{$\tilde{\alpha}$}
\psfrag{b}{$\beta$}
\includegraphics[scale = .5]{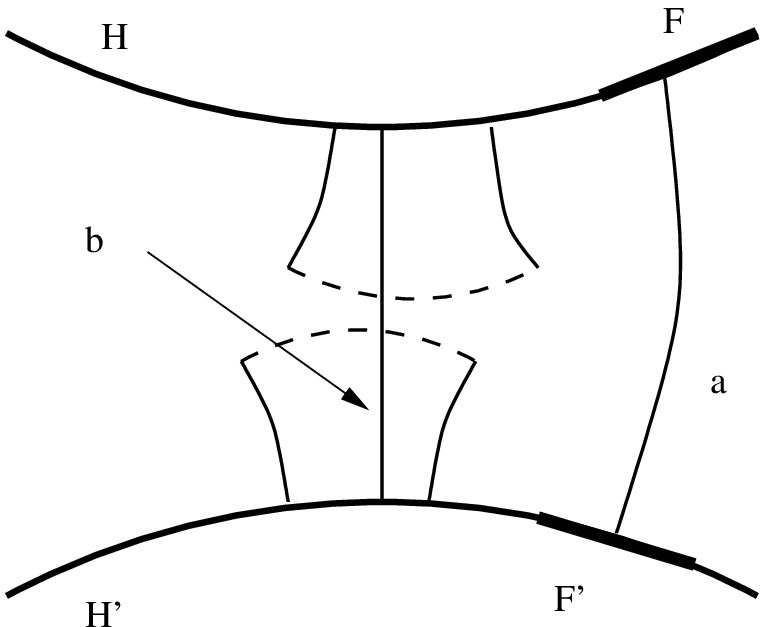}
\caption{\label{fig:ucover.fig}}
\end{center}
\end{figure}
Let $H$ and $H'$ be the hyperplanes spanned by $\tilde{F}$ and
$\tilde{F}'$, respectively.  Since $H$ and $H'$ are both mapped to the
face $F$ under the projection $\bbH^n\rightarrow\bbH^n/W=P$, they are
hyperparallel, thus have a (unique) common perpendicular $\beta$.
Moreover, $\beta$ must pass through at least two copies of the tile
$P$.  It follows from Lemma~\ref{lem:connectingpath} that the length
of $\beta$ is $\geq 2L$ with equality holding only if $\beta$ is 
the lift of an edge.  Since the geodesic $\tilde{\alpha}$ is at least
as long as $\beta$ with equality holding only if
$\alpha=\tilde{\beta}$, the closed geodesic $\alpha$ will have minimal
length only if it is an edge loop.

Now suppose $\phi:M_{\lambda}\rightarrow M_{\mu}$ is an isometry.  Since
$\phi$ takes minimal length closed geodesics to minimal length closed
goedesics, it must take edge loops to edge loops.  The 
$0$-cells in a small cover can be characterized as the points where
edge loops intersect, thus $\phi$ must take $0$-cells to $0$-cells
and, therefore, $1$-cells to $1$-cells.  Since any cell of
dimension $\geq 2$ in a small cover can be characterized as the convex
hull of its bounding $1$-cells, the isometry $\phi$ must take cells to
cells.  Thus, by Proposition~\ref{prop:equivalence}, $M_{\lambda}$ and
$M_{\mu}$ are equivalent.  
\end{proof}

Mostow rigidity gives the following:

\begin{corollary} 
If $P$ is the dodecahedron or the $120$-cell, then two small covers
of $P$ are homeomorphic if and only if they are equivalent. 
\end{corollary}

\section{Small covers of the dodecahedron and $120$-cell}\label{s:examples}

In this section, we describe an algorithm for enumerating equivalence
classes of small covers for general $P$, and apply it to the
dodecahedron and $120$-cell.  The complexity of the algorithm becomes 
unfeasible for the $120$-cell, so instead we apply it to a restricted
class of characteristic functions. 

\

\noindent{\bf The algorithm.}
Let $P$ be $n$-dimensional, and let $F_1,\ldots,F_d$ be any ordering
of the facets such that the first $n$ facets $F_1,\ldots,F_n$ all meet
at a vertex.  Given any characteristic function
$\lambda:W\rightarrow(\bbZ_2)^n$, we let $\lambda_i$ be the image
$\lambda(s_{F_i})$ of the fundamental reflection across the $i$th
facet.  By definition a $d$-tuple
$\blambda=(\lambda_1,\ldots,\lambda_d)$ determines a characteristic 
function if and only if the $\lambda_i$'s corresponding to facets
meeting at a vertex form a basis for $(\bbZ_2)^n$.  We call such a
$\blambda$ a {\em labeling of $P$} and each $\lambda_i$ a {\em
label}.  We say a labeling is {\em normalized} if the first $n$ labels
form the standard basis for $(\bbZ_2)^n$, and we let $\Lambda(P)$
denote the set of normalized labelings of $P$.  The following algorithm
determines $\Lambda(P)$.

\begin{algorithm}\label{alg:labelings}
Let $\lambda_1,\ldots,\lambda_n$ be the standard basis for
$(\bbZ_2)^n$, and let $S$ be the set $(\bbZ_2)^n-\{0\}$ (with some
fixed order).     
\begin{enumerate}
\item Assume by induction that $\lambda_1,\ldots,\lambda_{i}$ have
been chosen, and let $S_{i+1}$ be the list $S$ with all elements of
the following form removed: $\lambda_{i_1}+\ldots\lambda_{i_k}$
where that $i_j\leq i$ for $1\leq j\leq k$ and $F_{i_1}\cap\cdots\cap
F_{i_k}\neq\emptyset$.
\item If $S_{i+1}=\emptyset$ we know that there are no (unrecorded)
labelings that begin $\lambda_1,\ldots\lambda_{i}$.  In this case, we
back up until we find a nonempty $S_{l}$ with $l<i$ (and $l$ as large as 
possible).  We let $\lambda_l$ be the first element in $S_l$ and
go back to Step 1, using $l$ instead of $i$.  If
$S_{i+1}\neq\emptyset$, we let $\lambda_{i+1}$ be the first 
element in the list $S_{i+1}$ and remove it from $S_{i+1}$.  If
$i+1<d$, we go back to Step 1, using $i+1$ instead of $i$.  If
$i+1=d$, we add the full labeling $\blambda$ to $\Lambda(P)$ and
repeat Step 2.  
\end{enumerate}   
\end{algorithm}

The group of symmetries $A(P)$ acts on $\Lambda(P)$, and the
equivalence classes of small covers are in bijection with
$A(P)$-orbits.  To determine these orbits, we choose a set of
generators for $A(P)$ and form the graph whose vertices are
elements of $\Lambda(P)$ and whose edges join any pair of vertices
that differ by a generator of $A(P)$.  One then applies any standard
graph algorithm to determine the connected components of this graph.

\

\noindent{\bf The dodecahedron.} Let $P^3$ be the regular
dodecahedron.  The symmetry group $A(P^3)$ is the Coxeter group with
diagram $\bH_3$: 
\begin{figure}[h]
\begin{center}
\psfrag{3}{$3$}
\psfrag{5}{$5$}
\psfrag{a}{\scriptsize $a$}
\psfrag{b}{\scriptsize $b$}
\psfrag{c}{\scriptsize $c$}
\psfrag{d}{\scriptsize $d$}
\includegraphics[scale = .5]{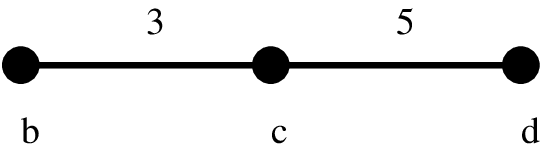}
\end{center}
\end{figure}

\begin{figure}[ht]
\begin{center}
\psfrag{1}{\scriptsize $1$}
\psfrag{2}{\scriptsize $2$}
\psfrag{3}{\scriptsize $3$}
\psfrag{4}{\scriptsize $4$}
\psfrag{5}{\scriptsize $5$}
\psfrag{6}{\scriptsize $6$}
\psfrag{7}{\scriptsize $7$}
\psfrag{8}{\scriptsize $8$}
\psfrag{9}{\scriptsize $9$}
\psfrag{10}{\scriptsize $10$}
\psfrag{11}{\scriptsize $11$}
\psfrag{12}{\scriptsize $12$}
\includegraphics[scale = .5]{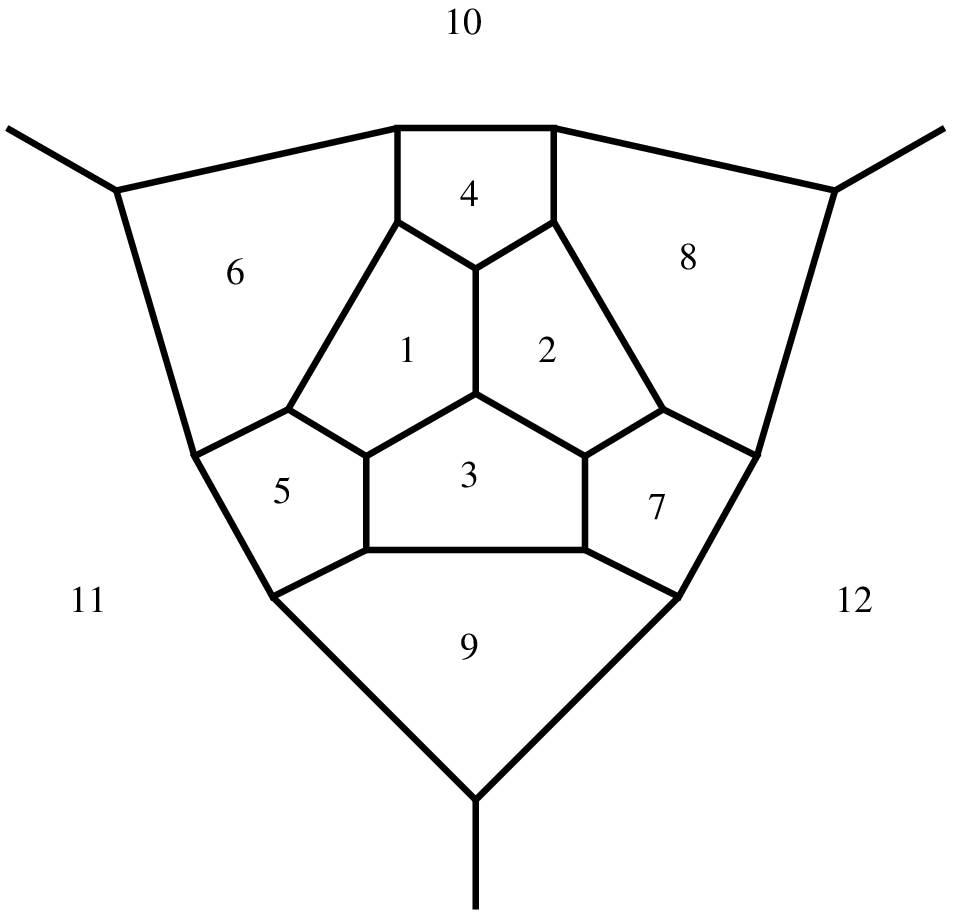}
\caption{\label{fig:dodec.fig}}
\end{center}
\end{figure}

\noindent Using the algorithm, we obtain $2165$ normalized labelings
of $P^3$.  Forming the graph on this set $\Lambda(P^3)$ with respect to the
standard generating set $\{b,c,d\}$ for $A(P^3)$, we obtain $25$
connected components.  Representative labelings are given in
Table~\ref{tab:dodec}.  To keep the data concise, we have used decimal
equivalents for $(\bbZ^2)^3$ (for example, $(1,0,0)\leftrightarrow 4$
and $(1,1,0)\leftrightarrow 6$).  The ordering we use for the facets
is shown in Figure~\ref{fig:dodec.fig}.  

\begin{center}
\begin{table}[ht]
\caption{Small covers of the dodecahedron \label{tab:dodec}}  
\begin{tabular}{||c|c||c|c||} \hline
\multicolumn{1}{||c|}{\raisebox{-2pt}{$\blambda=(\lambda_1,\ldots,\lambda_{12})$}} 
& \multicolumn{1}{|c||}{\raisebox{-2pt}{$A_{\lambda}$}}
&
\multicolumn{1}{|c|}{\raisebox{-2pt}{$\blambda=(\lambda_1,\ldots,\lambda_{12})$}} 
& \multicolumn{1}{|c||}{\raisebox{-2pt}{$A_{\lambda}$}}
\\[4pt]
\hline
\scriptsize$(1, 2, 4, 4, 2, 6, 1, 7, 7, 1, 3, 5)$ & $  1$ &
\scriptsize$(1, 2, 4, 4, 2, 6, 3, 5, 5, 3, 1, 7)$ & $  1$\\ \hline
\scriptsize$(1, 2, 4, 4, 2, 6, 3, 5, 5, 7, 3, 1)$ & $  1$ &
\scriptsize$(1, 2, 4, 4, 2, 6, 7, 1, 1, 3, 7, 5)$ & $  1$\\ \hline
\scriptsize$(1, 2, 4, 4, 2, 7, 1, 7, 7, 1, 3, 5)$ & $  1$ &
\scriptsize$(1, 2, 4, 4, 2, 7, 1, 7, 7, 1, 4, 2)$ &
$(\bbZ_2\times\bbZ_2)\ltimes\bbZ_3$\\ \hline 
\scriptsize$(1, 2, 4, 4, 2, 7, 3, 5, 5, 2, 6, 1)$ & $  1$ &
\scriptsize$(1, 2, 4, 4, 2, 7, 3, 5, 5, 6, 4, 7)$ & $  1$\\ \hline
\scriptsize$(1, 2, 4, 4, 2, 7, 3, 7, 5, 2, 6, 1)$ & $  1$ &
\scriptsize$(1, 2, 4, 4, 2, 7, 3, 7, 5, 5, 3, 1)$ & $  \bbZ_2$\\ \hline
\scriptsize$(1, 2, 4, 4, 2, 7, 3, 7, 5, 5, 6, 1)$ & $  1$ &
\scriptsize$(1, 2, 4, 4, 2, 7, 7, 1, 1, 2, 6, 5)$ & $  \bbZ_2$\\ \hline
\scriptsize$(1, 2, 4, 4, 2, 7, 7, 1, 5, 6, 3, 4)$ & $  1$ &
\scriptsize$(1, 2, 4, 4, 2, 7, 7, 3, 5, 5, 3, 1)$ & $  \bbZ_2\times\bbZ_2$\\ \hline
\scriptsize$(1, 2, 4, 4, 2, 7, 7, 3, 5, 5, 6, 1)$ & $  1$ &
\scriptsize$(1, 2, 4, 4, 3, 6, 3, 5, 5, 7, 2, 1)$ & $  1$\\ \hline
\scriptsize$(1, 2, 4, 4, 3, 6, 3, 7, 5, 5, 2, 1)$ & $  \bbZ_2$ &
\scriptsize$(1, 2, 4, 4, 3, 6, 3, 7, 6, 5, 2, 1)$ & $  1$\\ \hline
\scriptsize$(1, 2, 4, 4, 3, 6, 5, 3, 6, 5, 2, 1)$ & $  S_3$ &
\scriptsize$(1, 2, 4, 4, 3, 6, 7, 3, 2, 5, 7, 1)$ & $  \bbZ_2$\\ \hline
\scriptsize$(1, 2, 4, 4, 3, 6, 7, 3, 5, 1, 2, 6)$ & $  \bbZ_2$ &
\scriptsize$(1, 2, 4, 4, 3, 6, 7, 3, 5, 5, 2, 1)$ & $  \bbZ_2$\\ \hline
\scriptsize$(1, 2, 4, 4, 3, 7, 7, 3, 2, 5, 6, 1)$ & $  1$ &
\scriptsize$(1, 2, 4, 4, 3, 7, 7, 3, 5, 6, 2, 1)$ & $  \bbZ_2$\\ \hline
\scriptsize$(1, 2, 4, 5, 3, 7, 7, 3, 5, 4, 2, 1)$ &
$(\bbZ_2\times\bbZ_2)\ltimes\bbZ_6$ & & \\\hline 
\end{tabular} 
\end{table}
\end{center}

Also included in the table is the stabilizer subgroup
$A_{\lambda}\subset A(P^3)$ for the given labeling 
$\lambda$.  This subgroup acts via isometries on the corresponding
small cover $M_{\lambda}$ as does the group
$W/\Gamma_{\lambda}=(\bbZ_2)^3$.  It follows from Theorem~\ref{thm:main} 
that $M_{\lambda}$ admits no further symmetries.  In other words, the
isometry group of $M_{\lambda}$ is the semi-direct product of
$A_{\lambda}$ and $(\bbZ_2)^3$.  In particular, the small cover
on the bottom left has the largest symmetry group (with order $192$). 

\
  
\noindent{\bf The $120$-cell.} Let $P^4$ be the $120$-cell.  It has
$120$ dodecahedral facets, and the symmetry group $A(P^4)$ is the
Coxeter group with diagram $\bH_4$:  
\begin{figure}[h]
\begin{center}
\psfrag{3}{$3$}
\psfrag{5}{$5$}
\psfrag{a}{\scriptsize $a$}
\psfrag{b}{\scriptsize $b$}
\psfrag{c}{\scriptsize $c$}
\psfrag{d}{\scriptsize $d$}
\includegraphics[scale = .5]{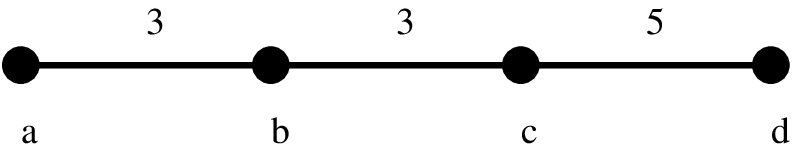}
\end{center}
\end{figure}

In theory, our algorithm can be used to search for all possible
labelings of $P^4$ with labels in the set $(\bbZ_2)^4-\{0\}$, but the
running time is too large to make the computation feasible.  It is not
hard to show that a labeling of $P^4$ must have at least $5$ labels,
and that if it has exactly $5$, they must be (modulo $G$) $1,2,4,8,15$
(decimal equivalents again).  Applying Algorithm~\ref{alg:labelings}
using $S=(1,2,4,8,15)$ (instead of $S=(\bbZ_2)^n-\{0\}$), the search
is effective and returns exactly $10$ labelings.  Another computation
shows that the stabilizer of one of these labelings has index $10$;
hence, all of these labelings are equivalent.  A representative
labeling $\blambda$ is given in Table~\ref{tab:120cell}.  

It follows from Theorem~\ref{thm:main} that up to homeomorphism
$M_{\lambda}$ is the unique small cover that uses only $5$ labels, and
that no small cover using more than $5$ labels is homeomorphic to
$M_{\lambda}$.  Since $A(P^4)$ has order $14,400$, the isometry group
of $M_{\lambda}$ has order $16\times 1,440=23,040$.     

\begin{center}
\begin{table}[ht]
\caption{\label{tab:120cell} A labeling for the
$120$-cell using $5$ labels.  The facets are indexed by their
barycenters, which were determined using the standard $4$-dimensional
geometric representation for $A(P^4)$.  This 
is the contragradient representation (\cite{B}, Chapter V, Section 4), 
so the generating involutions are linear reflections across the coordinate hyperplanes and are orthogonal with
respect to the {\em inverse} of the Coxeter matrix of type $\bH_4$.
With respect to this representation the vector $(2,0,0,0)$ is the
barycenter of the first facet and the remaining barycenters are the
translates of this point under the action of $A(P^4)$.}
\begin{tabular}{||c|c||c|c||} \hline
\multicolumn{1}{||c|}{\raisebox{-2pt}{center of facet $F$}} 
& \multicolumn{1}{|c||}{\raisebox{-2pt}{$\lambda(F)$}}
& \multicolumn{1}{|c|}{\raisebox{-2pt}{center of facet $F$}} 
& \multicolumn{1}{|c||}{\raisebox{-2pt}{$\lambda(F)$}}
\\
\hline
\scriptsize$\pm(2,0,0,0)$ & $ 1$ & \scriptsize$\pm(2,\sqrt{5}+1,-3-\sqrt{5},\sqrt{5}+1)$ & $ 1$\\ \hline
\scriptsize$\pm(-2,2,0,0)$ & $ 2$ & \scriptsize$\pm(-2,\sqrt{5}+3,-3-\sqrt{5},\sqrt{5}+1)$ & $ 8$\\ \hline
\scriptsize$\pm(0,-2,2,0)$ & $ 4$ & \scriptsize$\pm(2,\sqrt{5}+1,0,-1-\sqrt{5})$ & $ 4$\\ \hline
\scriptsize$\pm(0,0,-2,\sqrt{5}+1)$ & $ 8$ & \scriptsize$\pm(-2,\sqrt{5}+3,0,-1-\sqrt{5})$ & $ 8$\\ \hline
\scriptsize$\pm(0,0,\sqrt{5}+1,-1-\sqrt{5})$ & $ 8$ & \scriptsize$\pm(\sqrt{5}+3,-3-\sqrt{5},2,0)$ & $ 1$\\ \hline
\scriptsize$\pm(0,\sqrt{5}+1,-1-\sqrt{5},2)$ & $ 4$ & \scriptsize$\pm(\sqrt{5}+3,-1-\sqrt{5},-2,\sqrt{5}+1)$ & $ 8$\\ \hline
\scriptsize$\pm(0,\sqrt{5}+1,0,-2)$ & $ 15$ & \scriptsize$\pm(\sqrt{5}+1,-3-\sqrt{5},0,\sqrt{5}+1)$ & $ 4$\\ \hline
\scriptsize$\pm(\sqrt{5}+1,-1-\sqrt{5},0,2)$ & $ 2$ & \scriptsize$\pm(\sqrt{5}+3,-1-\sqrt{5},\sqrt{5}+1,-1-\sqrt{5})$ & $ 8$\\ \hline
\scriptsize$\pm(\sqrt{5}+1,-1-\sqrt{5},\sqrt{5}+1,-2)$ & $ 15$ & \scriptsize$\pm(\sqrt{5}+1,-3-\sqrt{5},\sqrt{5}+3,-1-\sqrt{5})$ & $ 2$\\ \hline
\scriptsize$\pm(\sqrt{5}+1,0,-1-\sqrt{5},\sqrt{5}+1)$ & $ 15$ & \scriptsize$\pm(\sqrt{5}+3,0,-1-\sqrt{5},2)$ & $ 2$\\ \hline
\scriptsize$\pm(\sqrt{5}+1,0,2,-1-\sqrt{5})$ & $ 2$ & \scriptsize$\pm(\sqrt{5}+1,0,-3-\sqrt{5},\sqrt{5}+3)$ & $ 4$\\ \hline
\scriptsize$\pm(\sqrt{5}+1,2,-2,0)$ & $ 8$ & \scriptsize$\pm(\sqrt{5}+3,0,0,-2)$ & $ 1$\\ \hline
\scriptsize$\pm(\sqrt{5}+3,-2,0,0)$ & $ 4$ & \scriptsize$\pm(\sqrt{5}+1,0,\sqrt{5}+1,-3-\sqrt{5})$ & $ 15$\\ \hline
\scriptsize$\pm(-1-\sqrt{5},0,0,2)$ & $ 1$ & \scriptsize$\pm(\sqrt{5}+1,\sqrt{5}+1,-1-\sqrt{5},0)$ & $ 15$\\ \hline
\scriptsize$\pm(-1-\sqrt{5},0,\sqrt{5}+1,-2)$ & $ 15$ & \scriptsize$\pm(2\sqrt{5}+2,-1-\sqrt{5},0,0)$ & $ 15$\\ \hline
\scriptsize$\pm(-1-\sqrt{5},\sqrt{5}+1,-1-\sqrt{5},\sqrt{5}+1)$ & $ 15$ & \scriptsize$\pm(-3-\sqrt{5},0,2,0)$ & $ 4$\\ \hline
\scriptsize$\pm(-1-\sqrt{5},\sqrt{5}+1,2,-1-\sqrt{5})$ & $ 4$ & \scriptsize$\pm(-3-\sqrt{5},2,-2,\sqrt{5}+1)$ & $ 2$\\ \hline
\scriptsize$\pm(-1-\sqrt{5},\sqrt{5}+3,-2,0)$ & $ 1$ & \scriptsize$\pm(-1-\sqrt{5},-2,0,\sqrt{5}+1)$ & $ 8$\\ \hline
\scriptsize$\pm(-3-\sqrt{5},\sqrt{5}+1,0,0)$ & $ 8$ & \scriptsize$\pm(-3-\sqrt{5},2,\sqrt{5}+1,-1-\sqrt{5})$ & $ 1$\\ \hline
\scriptsize$\pm(0,-1-\sqrt{5},0,\sqrt{5}+1)$ & $ 15$ & \scriptsize$\pm(-1-\sqrt{5},-2,\sqrt{5}+3,-1-\sqrt{5})$ & $ 8$\\ \hline
\scriptsize$\pm(0,-1-\sqrt{5},\sqrt{5}+3,-1-\sqrt{5})$ & $ 1$ & \scriptsize$\pm(-3-\sqrt{5},\sqrt{5}+3,-1-\sqrt{5},2)$ & $ 4$\\ \hline
\scriptsize$\pm(2,-3-\sqrt{5},\sqrt{5}+1,0)$ & $ 8$ & \scriptsize$\pm(-1-\sqrt{5},\sqrt{5}+1,-3-\sqrt{5},\sqrt{5}+3)$ & $ 1$\\ \hline
\scriptsize$\pm(-2,-1-\sqrt{5},\sqrt{5}+1,0)$ & $ 2$ & \scriptsize$\pm(-3-\sqrt{5},\sqrt{5}+3,0,-2)$ & $ 2$\\ \hline
\scriptsize$\pm(0,2,-3-\sqrt{5},\sqrt{5}+3)$ & $ 2$ & \scriptsize$\pm(-1-\sqrt{5},\sqrt{5}+1,\sqrt{5}+1,-3-\sqrt{5})$ & $ 15$\\ \hline
\scriptsize$\pm(2,-2,-1-\sqrt{5},\sqrt{5}+3)$ & $ 1$ & \scriptsize$\pm(-1-\sqrt{5},2\sqrt{5}+2,-1-\sqrt{5},0)$ & $ 15$\\ \hline
\scriptsize$\pm(-2,0,-1-\sqrt{5},\sqrt{5}+3)$ & $ 4$ & \scriptsize$\pm(0,-3-\sqrt{5},2,2)$ & $ 1$\\ \hline
\scriptsize$\pm(0,2,\sqrt{5}+1,-3-\sqrt{5})$ & $ 1$ & \scriptsize$\pm(0,-1-\sqrt{5},-2,\sqrt{5}+3)$ & $ 2$\\ \hline
\scriptsize$\pm(2,-2,\sqrt{5}+3,-3-\sqrt{5})$ & $ 4$ & \scriptsize$\pm(0,-3-\sqrt{5},\sqrt{5}+3,-2)$ & $ 4$\\ \hline
\scriptsize$\pm(-2,0,\sqrt{5}+3,-3-\sqrt{5})$ & $ 2$ & \scriptsize$\pm(0,-1-\sqrt{5},2\sqrt{5}+2,-3-\sqrt{5})$ & $ 15$\\ \hline
\scriptsize$\pm(0,\sqrt{5}+3,-1-\sqrt{5},0)$ & $ 2$ & \scriptsize$\pm(0,0,-3-\sqrt{5},2\sqrt{5}+2)$ & $ 8$\\ \hline
\end{tabular}
\end{table}
\end{center}

\end{document}